\documentclass[12pt]{amsart}
\usepackage{a4wide}

\usepackage{mathrsfs}
\usepackage{mathrsfs}
\def\xx#1 {\newtheorem{#1}[thm]{#1}}
\theoremstyle{definition}

\xx Lemma
\xx Theorem
\xx Definition
\xx Suggestion
\xx Construction
\xx Notation
\xx Remark
\xx Conclusion
\xx {Definition and Fact}
\xx Fact
\xx Assumption
\newcommand{\x}{{\tt x}}
\newcommand{\p}{{\tt p}}
\newcommand{\0}{{\tt 0}}

\title{Very many term clones\\ in a very  small variety}
\author{Martin Goldstern}
\address{DMG/Algebra, Technische Universit\"at Wien\endgraf
Wiedner Hauptstr 8-10\\1040 Wien, Austria (Europe)}
\email{Martin.Goldstern@tuwien.ac.at}

\date{2004-05-06; revised 2004-07-28}

\def\itm#1 {\item[{#1}]}

\begin{document}

\begin{abstract}
We give a nice  example of a finitely based locally finite variety 
which has uncountably many term clones. 
\end{abstract}

\maketitle

\newcommand{\M}{{\bf M}}
\newcommand{\E}{{\boldsymbol{ \epsilon}}}
\newcommand{\notE}{{\not\!\boldsymbol{ \epsilon}}}
\newcommand{\Ext}{{\bf Ext}}
\newcommand{\F}{{\mathfrak F}}
\newcommand{\A}{{\mathfrak A}}
\newcommand{\V}{{\mathbb V}}
\newcommand{\Ord}{{\bf Ord}}

\subsection*{Term clones}

Let $\A = (A, \Omega)$ be a universal  algebra. A {\em term function} is a
function $f:A ^n\to A$ (for some $n\in \{0,1,2,\ldots\}$) which is
induced by a term.   A {\em term clone of\/~$\A$} is a set of term functions 
which contains all the projections $\pi^n_k: A^n\to A$ and is closed
under composition  (also called ``superposition''). 
The {\em full term clone of\/~$\A$} is the set of all term functions.

Let $\V$ be a variety, $F_\V$ the free algebra in~$\V$ on countably 
many generators $\{x_1, x_2, \ldots\}$.
  A {\em term clone of\/~$\V$} is a term clone of~$F_\V$; 
since term functions are induced by elements of~$F_\V$, we can equivalently
view a term clone of~$\V$ as a   subset  $S$ of~$F_\V$ which 
contains all the generators and is closed under the following 
``substitution'' 
operation:
\begin{itemize}
\item[$(**)$]
 Whenever $t(x_1,\ldots, x_n)\in S$,
 and $t_1, \ldots, t_n\in S$, 
then also $t(t_1,  \ldots, t_n) \in S$. 
\end{itemize}

\subsection*{Our variety}
Ivan Chajda has asked whether there is a locally finite variety $\V$
 (preferably:
finitely based) which has uncountably many term clones.  We give here
a nice example\footnote{\'Agnes Szendrei
has pointed out that in fact many such varieties 
are known;  any minimal variety generated by a finite primal algebra (which 
has at least 3 elements) will have all the required properties.}
of such a variety.

Our language contains one binary operation symbol ${*}$ and two 
constant symbols $p$ and $0$.     We write $a{*} b {*} c$ for 
$(a{*} b) {*} c$.

  The laws of our variety $\V$ are
$$ 
   0{*} x=  x {*} 0=0,  \qquad 
       x{*} y{*} z=x{*} z{*} y,
\qquad 
     x{*} (y{*} z)=0,
\qquad 
          x{*} y{*} y  = 0.
$$  
%

\subsection*{The free algebra: elements} 
Let $\0,\p, \x_1,\x_2,\ldots$ be distinct objects. 

We will describe an algebra $\F = (F, {*}, p, 0)\in \V$  containing all $\x_i$ 
(in fact, $\F$ will be freely generated by the $\x_i$ in $\V$).

In addition to  $\0$, the set $F$  will contain the following distinct
objects: 
\begin{enumerate}
\item   Letters: $\p, \x_1, \x_2, \ldots$
\item  Words:     A word is a pair $w = (x,Y)$, where $ x$ is 
a letter and $Y$ is a finite nonempty set of letters.    Instead of
$w=(x,\{ y_1,\ldots, y_k\})$ (with all $y_i$ distinct),
we  write $w$ also as the string  $x\, y_1\, \cdots \,y_k$.  
We have to keep in mind that two apparently different strings
 such as~$x\,y_1y_2$ and $x\,y_2y_1$ are
only two different notations for the same
word $(x,\{y_1,y_2\})$.
\end{enumerate}
We define the {\em length} of $\0$ to be $0$, the length of any letter
is $1$, and the length of any word $(x,\{y_1,\ldots, y_k\})$ (with 
all $y_i$ distinct) is $k+1$. 

\subsection*{The free algebra: operations} 

\begin{itemize}
\item
The constant symbols $0$ and $p$ are interpreted as the objects $\0$ and
$\p$, respectively. 
\item
 The product is defined
naturally as follows: 
\begin{itemize}
\item If $x=0$ or $y=0$, then $x{*} y = 0$.
\item If $x$ and $y$ are letters, then $x{*} y = x\,y$.
\item If $y$ is a word, then $x{*} y = 0$. 
\item If $x$ is a word, say $x=x_0\,x_1\cdots  x_k$ with $k\ge 1$, and $y$ is 
a letter, $y\in \{x_1,\ldots, x_k\}$, then $x{*} y =0$. 
\item If $x$ is a word, say $x=x_0\,x_1\cdots  x_k$ with $k\ge 1$, and $y$ is 
a letter,  but $y\notin \{x_1,\ldots, x_k\}$, then $x{*} y = 
x_0\,x_1\cdots  x_k y$. 
\end{itemize}
\end{itemize}

It is easy to check (case by case) that 
  this operation yields an algebra in~$\V$, and it is also straightforward
to see that this algebra is free  over $\{\x_1, \x_2,\ldots\,\}$.

\subsection*{Local finiteness}
The free algebra over the empty set  has just 3 elements: $\{\0, \p, \p*\p\}$. 
The free algebra over one free generator  $\x$ is the set
$$
\{\0,\ \p,\ \x,\  \p\,\p,\ \p\,\x,\ \x\,\p,\  \x\,\x,\  \x\,\x\p=\x\,\p\x,\  \p\,\p\x=\p\,\x\p \}. $$
  In general, 
 the $\V$-free algebra over $n$ elements
has exactly $1 + (n+1)\cdot   2^{n+1}$ elements, so 
$\V$ is locally finite. 

\subsection*{Uncountably many term clones}

For any set  $ A \subseteq \{1,2,3,\ldots  \}$, we let 
  $S(A)$ consist of the element $\0$, plus the set of all
words that start with $p$ whose length is in $A$, i.e., 
all words of the form 
$$ \p \,\x_{i_2}\cdots \x_{i_n}  \qquad \mbox{ 
$i_2,\ldots, i_n$ are all distinct, and $n\in A$}
$$
or
$$ \p \, \p\x_{i_3}\cdots \x_{i_n}  \qquad \mbox{ 
$i_3,\ldots, i_n$ are all distinct, and $n\in A$}
$$
Applying the operation $(**)$ to any element $w \in S(A)$
will either result in~$0$, or in a word or letter $w'$ of the same length
as $w  $.  

{\small  [Applying $(**)$ to a word of the form $ w=x \,y_1\cdots y_n$, where 
$x$ is one of the letters $\x_i$,
may of course change the length of $w$; but the
first letter of any word in $S(A)$ is the constant letter $\p$, so all
applications of $(**)$ to a word in $S(A)$ will either just rename variables,
or produce $0$ because of $x*(y*z)=0$.]}

 So $S(A)$ is closed under the operation $(**)$. For $A\not=A'$ we 
have $S(A)\not= S(A')$.  So each $S(A)$ induces a different clone
of term functions.

\end{document}